\documentclass[10pt,twoside,final]{amsart}

\usepackage[english]{babel}
\usepackage{amsfonts,fancyhdr,graphics,amsmath,amssymb}

\title{Bireflectionality in the commutator subgroup of a finite orthogonal group}

\newtheorem{definition}{Definition}[section]
\newtheorem{theorem}{Theorem}[section]
\newtheorem{proposition}[theorem]{Proposition}
\newtheorem{lemma}[theorem]{Lemma}
\newtheorem{corollary}[theorem]{Corollary}
\newtheorem{remark}[theorem]{Remark}

\newcommand {\ch }{\operatorname{char}}
\newcommand {\SOG }{\operatorname{SO}}
\newcommand {\OG }{\mathrm{O}}

\newcommand {\GL }{\mathrm{GL}}
\newcommand {\SL }{\mathrm{SL}}
\newcommand {\End }{\mathrm{End}}

\newcommand {\Bahn }{\mathrm{B}}
\newcommand {\Fix }{\mathrm{Fix}}

\newcommand {\Neg }{\mathrm{Neg}}
\newcommand {\GF }{\mathrm{GF}}

\newcommand {\Idm }{\mathrm{I}}
\newcommand {\Cent }{\mathrm{Cent}}
\newcommand {\disc }{\operatorname{disc}}

\newcommand {\Wind }{\operatorname{ind}}

\usepackage{ifdraft}
\ifdraft{\usepackage[outer]{showlabels}
	\usepackage{todonotes} \usepackage{datetime}}{}
\usepackage[pagebackref,colorlinks,citecolor=red,urlcolor=blue,linkcolor=green, bookmarks=false,hypertexnames=true]{hyperref}
\usepackage{orcidlink}
\usepackage{fancyhdr}
\usepackage{verbatim}

\begin{document}

\bibliographystyle{plain}

\setcounter{page}{1}

\thispagestyle{empty}


\author{Klaus Nielsen}\,\orcidlink{0009-0002-7676-2944}
\email{klaus@nielsen-kiel.de}
\ifdraft{\today \ \currenttime}{\date{October 30, 2024}}

\pagestyle{fancy}
\fancyhf{}
\fancyhead[OC]{Klaus Nielsen}
\fancyhead[EC]{Bireflectionality in finite orthogonal groups}
%\fancyfoot[C]{\thepage}
\fancyhead[OR]{\thepage}
\fancyhead[EL]{\thepage}

\setlength{\headheight}{12.0pt}

\begin{abstract}
	Let $(V,f)$ be a finite-dimensional, nondegenerate, symmetric bilinear  space over a finite field $K = \GF(q)$ of characteristic not 2.
Let $\Omega(V)$ be the commutator subgroup of the orthogonal group $\OG(V) = \OG(V, f)$. It is shown that every element of $\Omega(V)$ is bireflectional (product of 2 involutions) if and only if it is reversible (conjugate to its inverse), except when   (i) $q \equiv 3 \mod 4$, (ii) $\dim V \equiv 2 \mod 4$, and $\disc V = -1$. Moreover, the bireflectional and reversible elements 
of $\Omega(V)$ are characterized.
\end{abstract}

\maketitle
%%%%%%%%%%%%%%%%%%%%%%%%%%%%%%%%%%%%%%%

%%%%%%%%%%%%%%%%%%%%%%%%%%%%%%%%%%%%%%%%%%%%%%%%%%%%%%%%%%%%%
\section{Introduction} \label{intro-sec}
Following H. Wiener\footnote{H. Wiener, Über die aus zwei Spiegelungen zusammengesetzten Verwandtschaften.
	Ber. Verh. kgl. Sächs. Ges. Wiss. Leipzig. Math.-phys. Cl. 43: 644-673 (1891)}, we call an element of a group $G$ bireflectional if it is a product of 2 involutions of $G$.
We say that an element of $G$ is reversible if it is conjugate to its inverse. Clearly, a bireflectional element is reversible. A group is bireflectional if all its elements are bireflectional.
 
We are interested in classifying reversible and bireflectional elements in subgroups of orthogonal groups. 
It is well known that the full orthogonal group $\OG(V) = \OG(V,f)$ is bireflectional by a theorem of Wonenburger \cite[Theorem 2]{Wonenburger-1966}. 
In \cite{KN-1987}, F. Knüppel and the author classified bireflectional elements in the special orthogonal group $\SOG(V)$:

\begin{theorem} \label{old-theorem}
	Let $V$ be a finite-dimensional, nondegenerate, symmetric bilinear space over an arbitrary field $K$ with characteristic  $\ch K \ne 2$.
	\begin{enumerate}
	 \item The special orthogonal group $\SOG(V)$ is bireflectional if and only if $\dim V \not \equiv 2 \mod 4$ or $V$ is a hyperbolic space over  $\GF(3)$.
	 \item Let $\dim V \equiv 2 \mod 4$ and $\varphi \in \SOG(V)$. Then $\varphi$ is bireflectional if and only if
	 $\varphi$ has an orthogonal summand of odd dimension.
\end{enumerate}
\end{theorem}

Later, they showed in \cite[Theorem 4.1]{KN-2010} that all reversible elements of $\SOG(V)$ are bireflectional.

In this paper, we consider the commutator subgroup $\Omega(V)$ of an orthogonal  group  over a finite field of odd characteristic and classify the bireflectional and reversible elements of $\Omega(V)$. Knüppel and Thomsen \cite[Theorem 8.5]{KT-1998} have already determined the bireflectional commutator groups:

\begin{theorem} \label{KT}
Let $K=\GF(q)$ be a finite field of odd characteristic. Then 
$\Omega(V)$ is bireflectional if and only if one of the following holds
\begin{enumerate}
	\item $q \equiv 1 \mod 4$ and  $\dim V \not \equiv 2 \mod 4$, or 
	\item  $\dim V \equiv 0 \mod 4$ and $\disc V = -1$, or 
	\item $\dim V \in \{8, 9\}$.
\end{enumerate}
\end{theorem}

Kim et al. \cite{KTV-2023} have shown that all reversibles of $\Omega(V)$ are bireflectional, except when $q \equiv 3 \mod 4$, $\dim V \equiv 2 \mod 4$, and  $\disc V = 1$.

Knüppel and the author \cite{KN-2010} also considered bireflectionality in
the commutator group of an orthogonal group over the reals. In a forthcoming paper, we classify bireflectional elements in the commutator subgroup of orthogonal groups of arbitrary Witt index over the reals.

%%%%%%%%%%%%%%%%%%%%%%%%%%%%%%%%%%%%%%%
\section{Main results}

We give a complete classification of reversible and bireflectional elements
of $\Omega(V)$. If $q \equiv 1 \mod 4$, then  all reversible elements of $\Omega(V)$ are bireflectional:

\begin{theorem} \label{main-th-1}
	Let $K=\GF(q)$ be a finite field with $q \equiv 1 \mod 4$.
	Let $\varphi \in \OG(V)$. Then the following are equivalent 
	\begin{enumerate}
		\item $\varphi^{\sigma} = \varphi^{-1}$ for some involution $\sigma \in \Omega(V)$;
		\item $\varphi^{\alpha} = \varphi^{-1}$ for some  element $\alpha \in \SOG(V)$;
		 \item $\varphi$ has an orthogonal summand of type $2\pm$ if $\dim V \equiv 2 \mod 4$.
	\end{enumerate}
\end{theorem}

We turn to the case the case $q \equiv 3 \mod 4$. First we classify the bireflectional elements of $\Omega(V)$.
We need a definition:
For $\varphi \in \End(V)$ let
$\Bahn^{\infty}(\varphi^2)$ be the Fitting one space of $\varphi^2-1$ and
$n_j(\varphi)$ denote the number of elementary divisors $(x\pm1)^d$ of $\varphi$ of degree $d \equiv j \mod 8$.

\begin{theorem} \label{main-th-2}
	Let $K = \GF(q)$ be a finite field with $q \equiv 3 \mod 4$.
	Let $\varphi \in \Omega(V)$. Then
	$\varphi$ is bireflectional in $\Omega(V)$ if and only if $\varphi$
	is bireflectional in $\SOG(V)$ and
	\begin{enumerate}
		\item $\varphi$ has an elementary divisor $p^d$, where $p \ne x\pm 1$ is irreducible, and $d$ is odd, or
		\item $\varphi$ has an orthogonal summand of even dimension and discriminant -1 or
		\item  $n_3(\varphi) + n_5(\varphi) + \frac{n_2(\varphi)}{2} + \frac{n_6(\varphi)}{2}
		\equiv \frac{\dim \Bahn^{\infty}(\varphi^2)}{4} \mod 2$.
	\end{enumerate}
\end{theorem}

\begin{corollary} 
	Let $K = \GF(q)$ be a finite field with $q \equiv 3 \mod 4$.
	Assume that one of the following holds
	\begin{enumerate}
		\item $\dim V \equiv 0 \mod 4$ and $\disc V = -1$, or
		\item $\dim V = 8$ or $\dim V = 9$.
	\end{enumerate}
	Then $\Omega(V)$ is bireflectional.
\end{corollary}

Next, we determine the reversible elements of $\Omega(V)$.
We have the following improvement of \cite[Theorem 3.6]{KTV-2023}:

\begin{theorem}  \label{main-th-3}
	Let $K=\GF(q)$ be a finite field with $q \equiv 3 \mod 4$.
	Let $\varphi \in \Omega(V)$. Then  $\varphi$ is reversible but not bireflectional if and only if $\varphi$ has an orthogonal decomposition 
	$\varphi = \varphi_O \perp \varphi_E$, where 
	\begin{enumerate}
		\item $\varphi_O$ is bicyclic with elementary divisors $e_1 = e_2 = (x\pm 1)^{2t+1}$, 
		\item  $\disc \varphi_O = 1$, $\dim \varphi_E \equiv 4 \mod 8$, and 
		\item $\varphi_E$ has no elementary divisor $p^d$, where $p$ is irreducible, and $d$ is odd.
	\end{enumerate}
\end{theorem}

Finally, we show

\begin{theorem}  \label{main-th-4}
Let $K = \GF(q)$ be a finite field with $q \equiv 3 \mod 4$.
Let $\varphi \in \Omega(V)$.  Assume that $\varphi$ is reversed by an element $\eta \in \OG(V)$ with $\eta^2 = -1$.\footnote{Clearly, $\Omega(V)$ contains no such element if $\dim V \equiv 2 \mod 4$ and 2 is a nonsquare.} 
Then $\varphi$ is bireflectional in $\Omega(V)$.
\end{theorem}

\begin{corollary} 
	Let $K=\GF(q)$ be a finite field with $q \equiv 3 \mod 4$.
	Let $\varphi \in \Omega(V)$ be reversible. 
	Then $\varphi$ is bireflectional in $\Omega(V)$ if and only its image in 
	$\mathrm{P}\Omega(V)$ is bireflectional.
\end{corollary}
%%%%%%%%%%%%%%%%%%%%%%%%%%%%%%%%%%%%%%%
\section{Preliminaries}
 
 %%%%%%%%%%%%%%%%%%%%%%%%%%%%%%%%%%%%%%
 \subsection{Notations}
 \mbox{}
 
 In this paper $(V,f)$ always denotes a finite-dimensional symmetric  bilinear space over a field $K$ of characteristic $\ne 2$. We assume that $f$ is nondegenerate.
 
 For a linear mapping $\varphi$ of $V$, let $\Bahn^j(\varphi)$ denote the image and $\Fix^j(\varphi)$ the kernel of $(\varphi -1)^j$. Put $\Fix^{\infty}(\varphi) = \bigcup_{j\ge 1} \Fix^j(\varphi)$ and $\Bahn^{\infty}(\varphi) = \bigcap_{j\ge 1}
 \Bahn^j(\varphi)$. The space $\Bahn(\varphi) := \Bahn^1(\varphi)$ is the path or residual space of $\varphi$.
 And $\Fix(\varphi) := \Fix^1(\varphi)$ is the fix space of $\varphi$. Further, let $\Neg^j(\varphi) =  \Fix^j(-\varphi)$. The space $\Neg(\varphi) := \Neg^1(\varphi)$ is the negative space of  $\varphi$.
 
 The reciprocal of a monic polynomial $r(x) \in K[x]$ of degree $d = \partial r$ is the polynomial $r^*(x) = r(0)^{-1}x^d r(x^{-1})$. A polynomial $r$ is selfreciprocal if $r = r^*$.
 
 By $\mu(\varphi)$ we denote the minimal polynomial of 
 $\varphi$. Clearly, the spaces $\Bahn^j(\varphi)$ and $\Fix^j(\varphi)$ are $\varphi$-invariant. A subspace  $W$ of $V$ is totally isotropic if $f(w, w) = 0$ for all $w \in W$. A subspace $W$ is totally isotropic if and only if it is totally degenerate, i.e. if $W \le  W^{\perp}$. The Witt index $\Wind V$ of $V$ is the dimension of a maximal totally isotropic subspace of $V$.
 By $\disc V$ we denote the discriminant of $V$.
 
 For a matrix $M \in K^{n,n}$ let $M'$ denote the transpose of $M$ and $M^+ = (M')^{-1}$ the transpose inverse of $M$.
 
 By $\Theta(\varphi)$ we denote the spinor norm  or spinorial norm of an orthogonal transformation $\varphi \in \OG(V)$, as defined by Lip\-schitz or Eichler. See Artin \cite[Definition 5.5]{EArtin}.
Other authors like Dieudonn\'e or Wall use a slightly different definition;
see e.g. \cite[ch. 9, 3.4 Definition]{WScharlau}. But when restricted to the special orthogonal group, all spinor norms are equal.

 We collect some well-known facts.
 
 \begin{remark} \label{remark-1}
 	We have $\det \varphi = \det \varphi|_{\Bahn(\varphi)}$.
 	If $\varphi \in \OG(V)$, then 
 	\begin{enumerate}
 		\item $\det \varphi = \det \varphi|_{\Neg(\varphi)} = \det \varphi|_{\Neg^{\infty}(\varphi)}$.
 		\item $\varphi$ is involutory if and only if $\Bahn(\varphi)  = \Neg(\varphi)$.
 		\item The invariant factors of $\varphi$ are selfreciprocal.
 		\item $Vr(\varphi)^{\perp} = \ker r^*(\varphi)$  for all $r \in K[x]$.
 		\item $V = \Bahn^{\infty}(\varphi) \oplus \Fix^{\infty}(\varphi)$.
 		\item If $\Bahn(\varphi)$ is totally isotropic, then $\dim \Bahn(\varphi)$ is even.
 		\item If $\sigma \in  \OG(V)$ is involutory, then $\Theta(\sigma) = \disc \Bahn(\sigma)$.
 		\item $\Omega(V) = \SOG(V) \cap \ker \Theta$ if $\Wind V \ge 1$.
 		\item If $\varphi$ is unipotent, then $\varphi \in \Omega(V)$.
 	\end{enumerate}	
 \end{remark}
 
 %%%%%%%%%%%%%%%%%%%%%%%%%%%%%%%%%%%%%%%
 \subsection{Orthogonally indecomposable transformations}
 \mbox{}
 
 Our proofs demand a closer look at orthogonally indecomposable orthogonal transformations. An orthogonal transformation  $\varphi \in \OG(V)$ (and also $V$) is called orthogonally indecomposable if $\varphi$ has no proper orthogonal summand; i.e. $V$ has no proper nondegenerate $\varphi$-invariant subspaces.
 
According to Huppert \cite[1.7 Satz]{Huppert-1980a}, an orthogonally indecomposable transformation of $\OG(V)$ is either\footnote{Our type enumeration follows Huppert \cite{Huppert-1980a}. In \cite{Huppert-1990} Huppert uses  a different  enumeration.}
\begin{enumerate}
	\item bicyclic with elementary divisors $e_1 = e_2 =(x \pm 1)^m$
	(type 1) or
	\item indecomposable as an element of $\GL(V)$ (type 2) or
	\item cyclic with minimal polynomial $(rr^*)^t$, where $r$ is irreducible and prime to its reciprocal $r^*$ ( type 3).
\end{enumerate}

It is convenient to introduce some subtypes for the types 2 and 3:

\begin{definition}
	Let $\varphi \in \OG(V)$ be orthogonally indecomposable of type 2 or 3 with 
	$\mu(\varphi) = p^m$, where $p$ is irreducible or $p = r r^*$, where $r$ is irreducible and prime to its reciprocal $r^*$.
	We say that $\varphi$ is of type 
	\begin{enumerate}
		\item $2-$ if $p = x - 1$,
		\item $2\pm$ if $p = x \pm  1$, 
		\item 2o if  $\varphi$ is of type 2, $p \ne  x \pm  1$, and $m$ is  odd,
		\item 2e if  $\varphi$ is of type 2, $p \ne  x \pm  1$, and  $m$ is even,
		\item 3o if $\varphi$ is of type 3, and $m$ is  odd,
		\item 3e if  $\varphi$ is of type 3, and $m$ is even.
	\end{enumerate}	
\end{definition}

In \cite[2.3 Satz]{Huppert-1980a} and \cite[3.2 Satz]{Huppert-1980b}, Huppert proved the following 

\begin{proposition} \label{prop-1}
	Let $\varphi \in \OG(V)$ be orthogonally indecomposable. Then 
	\begin{enumerate}
		\item If $\varphi$ is of type 1, then $\dim V \equiv 0 \mod 4$.
		\item $\dim V$ is odd if and only if $\varphi$ is of type $2\pm$.
	\end{enumerate}
\end{proposition}

\begin{definition}
Let $\varphi \in \OG(V)$. We say that $\varphi$ is
\begin{enumerate}
	\item of type T, if all orthogonally indecomposable  orthogonal  summands of $\varphi$ are of type T,
	\item  (type T)-free, if no orthogonally indecomposable  orthogonal  summand of $\varphi$ is of type T.
\end{enumerate}	
\end{definition}

\begin{lemma}      \label{lemma-disc}
	Let $\varphi \in \OG(V)$ be orthogonally indecomposable of type 2-.
	Assume that $\dim V \ge 3$. Then $\disc V = -\disc [\Bahn(\varphi)/ \Fix(\varphi)]$.
	\end{lemma}
	
	\begin{proof}
		Let $T$ be a complement of $\Fix(\varphi)$ in $\Bahn(\varphi)$. Then
		$T^{\perp}$ is a hyperbolic plane, and $T$ is isomorphic to $\Bahn(\varphi)/ \Fix(\varphi)$. Then $\disc V = (\disc T) (\disc T^{\perp}) = - \disc T = - \disc [\Bahn(\varphi)/ \Fix(\varphi)]$.
\end{proof}

%%%%%%%%%%%%%%%%%%%%%%%%%%%%%%%%%%%%%%%
\subsection{Products of 2 involutions}
\mbox{}

First of all, we outline a short proof of Wonenburger's theorems.
Let $P \in \GL(n, K)$ be reversible. Then $P$ is bireflectional: We may assume that $P$ is cyclic and that either $\Fix(P) =\Neg(P) = 0$  or $P$ is unipotent. Clearly, we may assume that $P$ is a Frobenius companion matrix. Then  $P$ is reversed by the antidiagonal unit matrix. 

Now let $\varphi \in \OG(V, f)$. To show that $\varphi$ is bireflectional, we may assume that $\varphi$ is orthogonally indecomposable. If $\varphi$ is cyclic, then any  involution $\sigma \in \GL(V)$ reversing $\varphi$ is an isometry of $(V,f)$, even if $f$ is degenerate:
Let $V = \langle u \rangle_{\varphi}$. Now $f(u \varphi^i \sigma, u\sigma \varphi^j \sigma ) = f(u \sigma\varphi^{-i}, u \varphi^{-j}) = f(u \sigma\varphi^j, u \varphi^i) = f(u \varphi^i, u\sigma \varphi^j)$.

So let $\varphi$ be of type 1.\footnote{ Wonenburger's proof of \cite[Lemma 5]{Wonenburger-1966} case contains an inaccuracy.} Then in a suitable basis

\[
\varphi = \left (\begin{array} {cc} A & 0\\ 0 & A^+ \end{array} \right ), 
f = \left (\begin{array} {cc} G & \Idm_m\\ \Idm_m & H\end{array} \right ).
\]
Let 
\[
\sigma = \left (\begin{array} {cc} S & 0\\ 0 & S' \end{array} \right ),
\]
where $A^S = A^{-1}$. Then $\varphi^{\sigma} = \varphi^{-1}$. As just shown $SGS' = G$ and $S'HS = H$ so that
$\sigma \in \OG(V)$. In fact, Huppert \cite[2.4 Satz]{Huppert-1980a} has shown that we even may assume that $G = H =0_m$. It follows that $\varphi = \pm \psi^2$ for some cyclic transformation $\psi \in \OG(V)$ with minimum polynomial $(x^2+1)^{2m}$. Clearly, if $\sigma$ reverses $\psi$, then $\sigma$ reverses $\varphi$.

For the convenience of the reader, we also provide a short proof of \cite[Theorem 4.1] {KN-2010}.

\begin{lemma}                                \label{lemma-cent}
	Let $\varphi \in \OG(V)$.
	Then the following are equivalent
	\begin{enumerate}
		\item $\Cent(\varphi) \subseteq \SOG(V)$;
		\item $\varphi$ is (type $2\pm$)-free.
	\end{enumerate}
\end{lemma}

\begin{proof}	
	Let $\xi \in \Cent(\varphi)$. Then $\Neg^{\infty}(\xi)$ is $\varphi$-invariant and nondegenerate. If $\det \xi = -1$, then $\dim \Neg^{\infty}(\xi)$ must be odd. 
	
	Conversely, let $T$ be a nondegenerate
	$\varphi$-invariant subspace of $V$. Let   $\sigma \in \OG(V)$ be the unique involution with $\Bahn(\sigma) = T$. Then  $\sigma$ commutes with $\varphi$. If $\dim T$ is odd, then $\det \sigma = -1$. 
\end{proof}

\begin{corollary}\label{cor-inv}
	Let $\varphi \in \SOG(V)$ be reversible. Then $\varphi$ is bireflectional.
\end{corollary}

\begin{proof}
	Let  $\varphi^{\alpha} = \varphi^{-1}$ for some $\alpha \in \SOG(V)$. 
	Let $\rho \in \OG(V)$ be an involution reversing $\varphi$. Assume that
	$\det \rho = -1$. Then $\xi := \alpha \rho$ commutes with $\varphi$, and 
	$\det \xi = -1$. By \ref{lemma-cent}, $\varphi$ has an orthogonal summand of type $2\pm$, and by \ref{old-theorem}, $\varphi$ is bireflectional.
\end{proof}

The following result reduces the problem of computing the determinant  and spinor norm of a reversing  involution to orthogonally indecomposable transformations:

\begin{lemma}                                \label{lemma-inv1}
	Let $\varphi \in \OG(V)$.  Assume that
	$\varphi$ is reversed by an involution $\sigma \in \OG(V)$.
	Then $V$ has a decomposition $V = U_1 \perp U_2 \perp \dots \perp U_m$, where the subspaces $U_j$ are $\varphi$- and $\sigma$-invariant and orthogonally indecomposable w.r.t. $\varphi$.
\end{lemma}

\begin{proof}
	See  \cite[Proposition, p. 212]{KN-1987a} or more general, \cite[Proposition 6.1]{GKN-2008}.
\end{proof}

First we deal with type $2\pm$ transformations. It suffices to consider the type 2-.
\begin{lemma}      \label{lemma-inv2}
	Let $\varphi \in \OG(V)$ be cyclic with minimal polynomial
	$(x-1)^{2m+1}$.  Let $\varphi = \sigma \tau$ be a product of
	two orthogonal involutions.  Then
	\begin{enumerate}
		\item $\dim \Bahn(\sigma) = \dim \Bahn(\tau) \in \{m,m+1\}$.
		\item Either $\Bahn(\sigma)$ and $\Bahn(\tau)$ or
		$\Fix(\sigma)$ and $\Fix(\tau)$ are hyperbolic spaces.
		\item $\Theta(\sigma) = \Theta(\tau) \in \{(-1)^t, (-1)^t \disc V \}$,
		where $t = \lfloor \frac{m+1}{2} \rfloor$.
	\end{enumerate}
\end{lemma}

\begin{proof}
	1: The fix space of $\varphi$ is invariant under $\sigma$ and $\tau$. Hence either 
	\begin{enumerate}
		\item $\Fix(\varphi) = \Fix(\sigma) \cap \Fix(\tau)$, $\Bahn(\sigma) \cap \Bahn(\tau) = 0$ or 
		\item $\Fix(\varphi) = \Bahn(\sigma) \cap \Bahn(\tau)$, $\Fix(\sigma) \cap \Fix(\tau) = 0$.
	 \end{enumerate}
	 as $\ch K \ne 2$. The assertion follows.
	
	2: The subspace $T = \Bahn^m(\varphi)$ is totally isotropic of dimension $m$ and
	$\langle \sigma, \tau \rangle$-invariant.
	Then  $ T = [T \cap  \Bahn(\sigma)] \oplus [T \cap  \Fix(\sigma)]$
	and  $ T = [T \cap  \Bahn(\tau)] \oplus [T \cap  \Fix(\tau)]$.
	We may assume that $\dim \Bahn(\sigma) = m = \dim \Bahn(\tau)$.
	If $m$ is even, then $\dim [T \cap  \Bahn(\tau)] = m/2 = \dim [T \cap  \Bahn(\sigma)]$, and $\Bahn(\tau)$ and  $\Bahn(\sigma)$ are hyperbolic.
	If $m$ is odd, then $\dim [T \cap \Fix(\tau)] = \frac{1}{2} (m+1) =
	\dim [T \cap \Fix(\sigma)]$, and $\Fix(\sigma)$ and $\Fix(\tau)$ are hyperbolic.
	
	Clearly, 3 follows immediately from 2.
\end{proof}

We turn to the hyperbolic types 1, 2e, and 3e.

\begin{lemma}  \label{lemma-inv3}
	Let $\varphi \in \OG(V)$ be orthogonally indecomposable of type 1, 2e or 3e.
	Let $\sigma \in \OG(V)$ be an involution reversing $\varphi$. Then 
	\begin{enumerate}
		\item $\Bahn(\sigma)$ and  $\Fix(\sigma)$ are hyperbolic spaces.
		\item $\dim \Bahn(\sigma) = \dim \Fix(\sigma) = \frac{1}{2} \dim V$.
		\item $\Theta(\sigma) = (-1)^{\frac{\dim V}{4}}$.\footnote{There is a typo in \cite[Lemma 5.5]{KT-1998}}
	\end{enumerate}
\end{lemma}

\begin{proof}
	1, 2: This is shown in  \cite[Lemma 2.4]{KhN-1987} and \cite[Lemma 2.9]{KhN-1987}. Clearly, 3 follows immediately from 1.
	See also \cite[Lemma 5.3]{KT-1998} and \cite[Lemma 5.5]{KT-1998}.
\end{proof}

Finally we consider the type 2o and 3o. Here we assume that $K$ is finite.

While the types 1, 2e, $2\pm$ and 3 have maximal Witt index, this is not always true in the type 2o case. Huppert \cite[4.1 Satz]{Huppert-1980b} has shown the following:

\begin{lemma} \label{lemma-inv4}
	Let $K$ be finite.
	Let $\varphi \in \OG(V)$ be orthogonally indecomposable of type 2o.
	Then $\Wind V = \frac{\dim V}{2} -1$.
\end{lemma}

We need 2 further  auxiliary results:

\begin{lemma}  \label{lemma-inv5}
	Let $K$ be finite.
	Let $\varphi \in \End(V)$. 
	Then the following are equivalent
	\begin{enumerate}
		\item All elementary divisors of $\varphi$ are squares;
		\item $\det \xi \in K^2$ for all $\xi \in \Cent(\varphi)$.
	\end{enumerate}
\end{lemma}

\begin{proof}
	Assume first that all elementary divisors of $\varphi$ are squares.
	Let $\xi \in \Cent(\varphi) \cap \GL(V)$.
	Considering the additive Jordan-Chevalley decomposition of $\varphi$ and the multiplicative  Jordan-Chevalley decomposition $\xi$, we may assume that $\varphi$ is nilpotent and $\xi$ is semisimple.
	Further, we may assume that $\xi$ is primary.
	Let $T = \bigoplus_{j=0}^n [\ker \varphi^j \cap  V \varphi^j]$.
	Then $T$ is obviously $\xi$-invariant, and $\dim T = \frac{\dim V}{2}$.
	Now $T$ has a $\xi$-invariant complement, as $\xi$ is semisimple. Hence
	$\det \xi \in K^2$.
	
	Conversely, assume  that $\varphi$ has an elementary divisor $p^t$.
	Let $\lambda \in K$. We may assume that $\varphi$ is cyclic. 
	Consider $\varphi$ in Jacobson normal form $M := \mathrm{I}_t \otimes P + \mathrm{J}_t \otimes \mathrm{N}$, where $P \in \GL(m, K)$ is irreducible with minimal polynomial $p$ and $N  \in K^{m,m}$ has rank one.
	It follows from  Newman \cite[Lemma 1]{Newman-1982} or Waterhouse \cite[Theorem 5]{Waterhouse-1984} or \cite[ch 2, 15.4 Theorem]{WScharlau} that there exists a polynomial $f$ such that  $\det f(P) = \lambda$. Then  $f(M) $ commutes with $M$ and $\det f(M) = \lambda^t$.
\end{proof}

	The next lemma is the matrix version of the following elementary fact:
If $\varphi \in \End(V)$ is selfadjoint w.r.t. a bilinear form $g$ on $V$, then 
$g$ is symmetric if $\varphi$ is cyclic.

\begin{lemma}  \label{lemma-inv5a}
	Let $P, G \in K^{n,n}$ with $PG = GP'$. If $P$ is cyclic, then $G$ is symmetric.
\end{lemma}

\begin{proof}
In fact, this is a well-known result.
 Frobenius has  shown that $\Cent(A) = K[A]$ and 
$A^S = A'$ for some nonsingular symmetric matrix $S$. Then $G = p(A) S$ for some polynomial $p$ so that $G' = S p(A') = p(A) S = G$.
\end{proof}

\begin{remark} \label{remark-2}
Let $\varphi = \rho \sigma$ be the product of involutions $\rho, \sigma \in \GL(V)$. If $\Fix(\varphi) = \Neg(\varphi) = 0$, then 
 $\dim \Bahn(\rho) = \dim \Fix(\rho) = \dim \Bahn(\sigma) = \dim \Fix(\sigma) = \frac{\dim V}{2}$.
\end{remark}

\begin{lemma} \label{lemma-inv6}
	Let $K$ be finite.
	Let $\varphi \in \OG(V)$ be orthogonally indecomposable of type 2o or 3o.
	\begin{enumerate}
		\item $\dim \Bahn(\rho) = \dim \Fix(\rho) = \frac{\dim V}{2}$ for all involutions $\rho \in \OG(V)$ inerting $\varphi$.
		\item There exists an involution  $\sigma \in \OG(V)$ reversing $\varphi$ with presrcibed spinor norm.
	\end{enumerate}
\end{lemma}

\begin{proof}
	1: see \ref{remark-2}
	
	2: If $\dim V \equiv 2 \mod 4$ apply  \cite[Lemma 5.8]{KT-1998} and \cite[Lemma 5.7]{KT-1998}. So let $\dim V \equiv 2 \mod 4$. If $\varphi$ is of type 2o, then by \ref{lemma-inv4}, $\Theta(\sigma) \ne \Theta(-\sigma)$ for all involutions $\sigma \in \OG(V)$ reversing $\varphi$. We are left with the case 3o. Then $\varphi$ is similar to a matrix $A \oplus A^+$, where $A$ and $A^+$ are coprime. It follows from \ref{lemma-inv5} that $A$ is conjugate to its transpose $A'$ via a matrix $S$ with  prescribed determinant.  
	By \ref{lemma-inv5a}, $S$ must be symmetric.

	In a suitable basis, we have
	\[
	\varphi = \left (\begin{array} {cc} A & 0\\ 0 & A^+ \end{array} \right ), 
		f = \left (\begin{array} {cc} 0 & \Idm\\ \Idm & 0\end{array} \right ).
	\]
	Let 
	\[
	\sigma = \left (\begin{array} {cc} 0 & S\\ S^{-1} & 0\end{array} \right ).
	\]
	
	Then $\sigma \in \OG(V)$, $\varphi^{\sigma} = \varphi^{-1}$, and $\Theta(\sigma) = \det S$.
\end{proof}

%%%%%%%%%%%%%%%%%%%%%%%%%%%%%%%%%%%%%%%
\section{Proofs}

We always assume that $K = \GF(q)$ is finite.

%**********************************************************************
%**********************************************************************
\subsection{Proof of theorem \ref{main-th-1}}
\mbox{}

\begin{proof}[Proof of theorem \ref{main-th-1}]

The implication $1 \longrightarrow 2$ is clear. 

Assume that 2 holds, and 
$\varphi$ is (type $2\pm$)-free. Then $\varphi \in \SOG(V)$.
By \ref{cor-inv}, $\varphi$ is bireflectional in $\SOG(V)$.
And by theorem \ref{old-theorem}, $\dim V \not \equiv 2 \mod 4$.

We prove  $3 \longrightarrow 1$: 

By \ref{old-theorem}, $\varphi$ is reversed by an involution $\sigma \in \SOG(V)$.
By  \ref{lemma-inv1} and \ref{lemma-inv6}, we can adjust the spinor norm without changing its determinant if $\varphi$ has an orthogonal summand of type 2o or 3o.

So let  $\varphi$ be (type 2o)- and (type 3o)-free.
We may assume that $\varphi$ is orthogonally indecomposable. By \ref{lemma-inv1}, \ref{lemma-inv2}, and \ref{lemma-inv3},
$\varphi$ is reversed by an involution $\tau \in \Omega(V)$.
\end{proof}
%**********************************************************************
%**********************************************************************
\subsection{Proof of theorem \ref{main-th-2}}
\mbox{}

We say that $\varphi$ is weakly homodisc if all orthogonally indecomposable summands in an orthogonal decomposition of $\varphi$ have the same discriminant.

\begin{lemma} \label{lemma-p1}
	Let $\varphi \in \Omega(V)$ be of type 2- and weakly homodisc.
	Let $\sigma \in \SOG(V)$ be an involution inverting $\varphi$. Then 
	$\Theta(\sigma) = (-1)^m$, where $m = n_3(\varphi) + n_5(\varphi)$.
\end{lemma}

\begin{proof}
		By induction, we may assume that $\varphi$ is either cyclic or bicyclic. Apply \ref{lemma-inv2}(3).
\end{proof}
	
	\begin{remark} \label{remark-p1}
		Let $\varphi \in \OG(V)$. Assume that $\varphi = \varphi_1 \perp \varphi_2$, where $\varphi_1$ and $\varphi_2$ are orthogonally indecomposable of type $2\pm$. Let $\disc V = -1$. Then there exists an involution $\sigma$ inverting $\varphi$ with prescribed determinant and prescribed spinor norm.
	\end{remark}
	
	\begin{proof}
		This is clear.
	\end{proof}
	
	\begin{proof}[Proof of theorem \ref{main-th-2}]
	 Let $\varphi \in \Omega(V)$ be bireflectional in $\SOG(V)$.
	 
	 First  assume that $\varphi$ is not (type 2o)- or (type 3o)-free. 
	 Then we can adjust the spinor norm of any involution $\sigma \in \SOG(V)$ reversing $\varphi$ without changing its determinant.
	 
	 Next assume that $\varphi$ has an orthogonal summand $\psi$ with $\dim  \psi \equiv 0 \mod 2$ and $\disc  \psi = -1$.
	 We may assume that $\varphi$ is (type 2o)- and (type 3o)-free.
	 Then $\psi = \psi_1 \perp \psi_2$, where  $\psi_1$ is orthogonally indecomposable of type $2\pm$ , and $\disc \psi_1 = -1$, as $\disc  \psi = -1$.
	 Further $\psi_2$ must have an  orthogonally indecomposable orthogonal summand of type $2\pm$ with discriminant 1, as $\dim \psi$ is even.
	 Now it follows from \ref{remark-p1} that $\varphi$ is bireflectional.
	 
So assume that (1) and (2) do not hold.
Let $\sigma \in \SOG(V)$ be an involution  reversing $\varphi$.
By \ref{lemma-inv1}, there exists a decomposition 
$V = E \perp O$, where $O$ and $E$ are invariant under $\varphi$ and $\sigma$, 
$E$ is (type $2\pm$)-free, and every orthogonally indecomposable  summand of $O$ is of type $2\pm$.
By \ref{lemma-p1}, $\Theta(\sigma_O) = (-1)^m$, where $m = n_3(\varphi) + n_5(\varphi)$. By \ref{lemma-inv3}, $\Theta(\sigma_E) = (-1)^d$, where $d = \frac{\dim E}{4} = \frac{n_2(\varphi)}{2} + \frac{n_6(\varphi)}{2}+ \frac{\dim \Bahn^{\infty}(\varphi^2)}{4}$.
\end{proof}

%**********************************************************************
%**********************************************************************
\subsection{Proof of theorem \ref{main-th-3}}
\mbox{}

We call a linear transformation $\varphi$ heterocyclic if all of its elementary divisors have multiplicity one. 
We say that an
orthogonal transformation $\varphi \in \OG(V)$ is homodisc if all of its orthogonally indecomposable orthogonal summands have the same discriminant.

\begin{remark}
	Let $\varphi \in \OG(V)$ be unipotent of type 2-. Assume that $\varphi$ has an elementary divisor $(x-1)^{2m+1}$ of muliplicity one. Then
	all $\varphi$-invariant subspaces $U$ of $V$ of dimension $2m+1$ have discriminant $(-1)^m \dim [\Fix^{m+1}(\varphi) \cap \Bahn^{m}(\varphi)]/[\Fix^{m}(\varphi)  + \Bahn^{m+1}(\varphi)]$ by \ref{lemma-disc}.
\end{remark}

\begin{corollary} 
	Let $\varphi \in \OG(V)$ be of type $2\pm$ and heterocylic. Then
	$\varphi$ has an orthogonal summand of even dimension and discriminant -1 if and only if $\varphi$ is not homodisc. 
\end{corollary}

\begin{lemma}                           \label{lemma-p3}
	Let $\varphi \in \OG(V)$ be of type $2\pm$. If $\varphi$ is not heterocyclic, then  $V$ has a nondegenerate $\varphi$-invariant subspace $U$ with prescribed discriminant.
\end{lemma}

\begin{proof}
	We may assume that $\varphi$ has elementary divisors $e_1 = e_2 = (x-1)^{2t+1}$. Then
	$\Bahn^t(\varphi)/\Fix^t(\varphi)$ is
	a nondegenerate plane $H$. So there exists a vector $b \in \Bahn^t(\varphi)$ with prescribed $f(b,b)$ as $H$ is universel. There exists a vector  $u \in V$ such that $u(\varphi-1)^t = b$. Then $U = \langle u\rangle_{\varphi}$ is nondegenerate and $\disc U = (-1)^t f(b,b)$
	by  \ref{lemma-disc}.
\end{proof}

\begin{lemma}                                \label{lemma-p3a}
 Let $\varphi \in \OG(V)$. Then the following are equivalent
	\begin{enumerate}
		\item $\Cent(\varphi) \subseteq \Omega(V)$;
		\item $\varphi$ is (type $2\pm$)-, (type 2o)-, and (type 3o)-free.
	\end{enumerate}
\end{lemma}

\begin{proof}
	If $\varphi$ is not (type 2o)- or (type 3o)-free, then 
	$\Cent(\varphi) \nsubseteq \ker \Theta$ by \ref{lemma-inv6}.
	If $\varphi$ is not (type $2\pm$)-free, then $\Cent(\varphi) \nsubseteq \SOG(V)$ by \ref{lemma-cent}. 
	
	Conversely, let $\varphi$ be (type $2\pm$)-, (type 2o)-, and (type 3o)-free.
	We already know that $\Cent(\varphi) \subseteq \SOG(V)$ by \ref{lemma-inv6}.
	It remains to show that $\Cent(\varphi) \subseteq \ker \Theta$. Let $\xi \in \Cent(\varphi)$. Considering the Jordan-Chevalley decompositions of $\varphi$ and $\xi$, we may assume that $\xi$ is semisimple and $\varphi$ is unipotent. We may further assume that $\xi$ is primary. Let $V = U_1 \perp \dots \perp U_m$ be the orthogonal summand of (orthogonally) indecomposable $\xi$-invariant subspaces.
	We show that $m$ must be even. First let $\xi$ be of type 2. By \ref{lemma-inv4}, the subspaces 
	$U_j$ are not hyperbolic. Hence $m$ is even, as $V$ is hyperbolic.
	So let $\xi$ be of type 3o with minimal polynomial $rr^*$, where
	$r$ is  prime to its reciprocal $r^*$. Let $T = \ker r(\xi)$.
	Then $T$ is $\varphi$-invariant. Let $\psi = (\varphi-1)|_T$. Let  $S = [T \psi \cap \ker \psi] + [T \psi^2 \cap \ker \psi^2] + \dots$. Then $S$
	is $\xi$-invariant, and $\dim S = \frac{\dim T}{2}$ as all elementary divisors of $\psi$ have even degree. Again $m$ must be even. It follows that $\Theta(\xi) = 1$.
\end{proof}

\begin{lemma}                                \label{lemma-p3b}
	Let $\varphi = \varphi_O \perp \varphi_E \in \OG(V)$, where
	$\varphi_O$ is of type $2\pm$ and $\varphi_E$ is of type 1. Assume further that $\varphi_O$ is homodisc with common discriminant $\delta$. 
	\begin{enumerate}
		\item If $\delta = 1$, then $\Cent(\varphi) \subseteq \ker \Theta$.
		\item If $\delta = -1$, then $\Cent(\varphi) \cap \SOG(V) \subseteq \Omega(V)$.
	\end{enumerate}
\end{lemma}

\begin{proof}
	Let $\xi \in \Cent(\varphi)$. Let $V = F \perp N \perp U$, where $F = \Fix^{\infty}(\xi)$ and $N = \Neg^{\infty}(\xi)$. Clearly, $F, N$ and $U$ are $\varphi$-invariant, and $\Theta(\xi_F) = 1$. Clearly, $\Theta(\xi_N) = \Theta(- \xi_N) \disc N = 1$ if $\delta = 1$. If $\delta = -1$, then  $\Theta(\xi_N) = \Theta(- \xi_N) \disc N = \det \xi_N = 1$. 
	
	Suppose that $\Theta(\xi) = -1$. Then $U \ne 0$. So we may assume that $V = U$.
	By \ref{lemma-p3a}, $\varphi \ne \varphi_E$, and by \ref{lemma-p3}, $\varphi_O$ contains exactly one elementary divisor $e = (x-1)^{2t+1}$ of maximal odd degree. Then $\Fix^{t+1}(\varphi)/\Fix^t(\varphi)$
	is  1-dimensional, nondegenerate,  and $\xi$-invariant, a contradiction.
\end{proof}

\begin{lemma}                                \label{lemma-p3c}
Let  $\disc V = 1$. Let $\varphi \in \OG(V)$
	be  bicyclic with  elementary divisors $e_1 = e_2 = (x \pm 1)^{2m+1}$.  There exist involutions $\sigma_1, \sigma_2, \sigma_3 \in \OG(V)$ and an element
	$\psi \in \OG(V)$
	reversing $\varphi$ such that
	\begin{enumerate}
		\item  $\Theta(\sigma_1) = 1, \det \sigma_1 = 1$,
		\item  $\Theta(\sigma_2) = 1, \det \sigma_2 = -1$,
		\item $\Theta(\sigma_3) = -1,   \det \sigma_3 = -1$,
		\item $\Theta(\psi) = -1,   \det \psi = 1$.
	\end{enumerate}
\end{lemma}

\begin{proof}
	By \ref{lemma-p3}, there exists a decomposition $V = U \perp W$, where $U$ and $W$ are $\varphi$-cyclic, and $\disc U$ can be prescribed. If $\disc U = 1$, then by \ref{lemma-inv2},
	there exist involutions $\sigma_1, \sigma_2 \in \OG(V)$ reversing $\varphi$ with $\Theta(\sigma_1) = \det \sigma_1 = 1$
	and $\Theta(\sigma_2) = -\det \sigma_2 = 1$.
	If $\disc U = -1$, then again by \ref{lemma-inv2},
	there exists an involution $\sigma_3 \in \OG(V)$ reversing $\varphi$ with $\Theta(\sigma_3) = \det \sigma_3 = -1$.
	Put $\psi = \sigma_1 \sigma_2 \sigma_3$.
\end{proof}

\begin{proof}[Proof of theorem \ref{main-th-3}]
	Let $\varphi \in \Omega(V)$ be reversible, and let $\alpha \in \Omega(V)$ be an element reversing $\varphi$.
	Suppose that $\varphi$ is not bireflectional.
	It follows from \ref{main-th-2} that $\varphi$ must be 2o-free and 3o-free.
	By \ref{cor-inv}, $\varphi$ is reversed by an involution $\sigma \in \SOG(V)$.  Then $\Theta(\sigma) = -1$.
	
	Let $\varphi = \varphi_O \perp \varphi_E$ be any orthogonal decomposition, where $\varphi_O$ is of type $2\pm$ and $\varphi_E$ is (type $2\pm$)-free.
	By \ref{lemma-p3a}, $\varphi \ne \varphi_E$ : otherwise $\alpha \sigma \in \Cent(\varphi) \cap \SOG(V) \subseteq \ker \Theta$. Further, by  \ref{lemma-p3b}, $\varphi_O$ is not homodisc. Hence $\varphi_O$ must have at least 2 elementary divisors.
	By \ref{main-th-2}, $\varphi_O$ must be weakly homodisc. Thus 
	$\varphi_O$ is not  heterocyclic. By \ref{lemma-p3}, $\varphi_O$ must be bicyclic.
	By \ref{main-th-2}, $\varphi$ is not bireflectional if and only if
	$\dim \varphi_E \equiv 4 \mod 8$.
	
	Conversely, it follows from \ref{lemma-p3c} that $\varphi$ is reversible.
\end{proof}

%**********************************************************************
%**********************************************************************
\subsection{Proof of theorem \ref{main-th-4}}

\begin{lemma}                                \label{lemma-p4}
	Let $\varphi \in \OG(V)$. Let $\varphi$ be (type 1)-free if $-1 \in K^2$. Assume that
	$\varphi$ is reversed by an element   $\eta \in \OG(V)$ with $\eta^2 = -1$.
	Then $V$ has a decomposition $V = U_1 \perp U_1\eta \perp \dots \perp U_m \perp U_m \eta$, where the subspaces $U_j$ are orthogonally indecomposable $\varphi$-invariant subspaces.
\end{lemma}

\begin{proof}
	Let $T$ be a  not necessarily nondegenerate $\varphi$-cyclic subspace of $V$. Then $T \eta \le T^{\perp}$.
	
	So by induction, we may assume that  $\varphi$ is of type 1 (and unipotent). 
	Let $u \in V - \Fix^{t-1}(\varphi)$, where  $t = \partial \mu(\varphi)$. Put $U = \langle u \rangle_{\varphi}$. Then $U \cap U \eta = 0$ as $-1 \not \in K^2$. Let $w \in (U \eta)^{\perp}- [u(\varphi -1)^{t-1}]^{\perp}$. Put $W = \langle w \rangle_{\varphi}$. Then $U\oplus W$ is orthogonally indecomposable and $U\eta  \oplus W\eta \le [U\oplus W]^{\perp}$. We are done by induction.
\end{proof}

\begin{proof}[Proof of theorem \ref{main-th-4}]
	 If $\varphi$ is not (type 2o)- or not (type 3o)-free, then $\varphi$ is bireflectional by \ref{main-th-2}(1). So let $\varphi$ be (type 2o)- and (type 3o)-free. By \ref{lemma-p4}, $\dim \Bahn^{\infty}(\varphi^2) \equiv 0 \mod 8$, and $n_3(\varphi), n_5(\varphi), 
	 \frac{n_2(\varphi)}{2}$, and $\frac{n_6(\varphi)}{2}$ are even.
	By \ref{main-th-2}(3), $\varphi$ is bireflectional.
\end{proof}

%%%%%%%%%%%%%%%%%%%%%%%%%%%%%%%%%%%%%%%
\section{Appendix. Bireflectionality in  $\SL(V)$}
\mbox{}

We want to prove the following result:

\begin{theorem} \label{app-th1}
	Let $\ch K \ne 2$.
	Let $\varphi \in \SL(V)$ be similar to its inverse in $\GL(V)$. 
	Then the following are equivalent
	\begin{enumerate}
		\item  $\varphi$ has an elementary divisor $(x \pm 1)^d$ of odd degree $d$ if $\dim V \equiv 2 \mod 4$;
		\item $\varphi$ is bireflectional in $\SL(V)$.
	\end{enumerate}
\end{theorem}

If $K$ is finite, this has already been proven by Gill and Singh
\cite[Theorem 6.1]{GillSingh-2011}. 

For the proof of \ref{app-th1}, we use the following 
unpublished result of Thomsen \cite[1.9 Satz]{Thomsen1991}:

\begin{proposition} \label{prop-2}       
	Let $\ch K \ne 2$. Let $\varphi \in \GL(W)$, and let $\sigma \in
	\GL(W)$ be an involution with $\varphi^{\sigma} = \varphi^{-1}$.
	Then $W = W_1 \oplus \dots \oplus W_m$,
	where the subspaces $W_j$ are $\sigma$-invariant and $\varphi$-invariant, and (i) $W_j$ is $\varphi$-indecomposable or (ii) $W_j = X_j \oplus X_j \sigma$, where $X_j$ is $\varphi$-indecomposable.
\end{proposition}

We provide a short proof. First we show

\begin{lemma} \label{lemma-34a}                             
	Let $\ch K \ne 2$. Let $\varphi \in \GL(W)$ be primary,
	and let $\sigma \in
	\GL(W)$ be an involution with $\varphi^{\sigma} = \varphi^{-1}$.
	Let  $U$ be a $\varphi$-cyclic subspace of maximum dimension.  Assume that $U$ is $\sigma$-invariant. Then
	$U$ has a $\langle \varphi, \sigma \rangle$-invariant complement.
\end{lemma}

\begin{proof}
	Let $\mu_{\varphi} = p^s$, where $p$ is irreducible.
	We may assume that $\varphi$ is not cyclic. Then
	$\varphi$ is not cyclic on $\ker p(\varphi)$.
	Furthermore, $\ker p(\varphi)$ is $\sigma$-invariant.
	Hence there exists a minimal $\langle \varphi, \sigma \rangle$-invariant subspace  $M$ such that $U \cap M = 0$. Then by induction
	$V/M = [U\oplus M]/M \oplus C/M$ for some $\langle \varphi, \sigma \rangle$-invariant subspace $C$.
	Hence $V = U \oplus C$, and we are done.
\end{proof}

\begin{proof}[Proof of \ref{prop-2}]
	We may assume that $\mu(\varphi) = p^s$ or $\mu(\varphi) = (pp^*)^t$ , where $p$ is irreducible, and  $p$ is either selfreciprocal or prime to its reciprocal $p^*$. 
	The biprimary case is easy: 
	Let $\ker p^t = U_1 \oplus \dots \oplus U_m$, where the subspaces $U_j$ are $\varphi$-cyclic. Then $\ker (p^*)^t = U_1 \sigma \oplus \dots \oplus U_m \sigma$,
	and $W_j = U_j \oplus U_j \sigma$ is cyclic and $\sigma$-invariant.
	So let $\varphi$ be primary. Let  $v \in [\Fix(\sigma) \cup \Neg(\sigma)] - \ker p^{t-1}(\varphi)$. Then $\langle v \rangle_{\varphi}$ is $\sigma$-invariant.
	By \ref{lemma-34a}, $\langle v \rangle_{\varphi}$ has a $\langle \varphi, \sigma \rangle$-invariant complement. We are done by induction.
\end{proof}	

\begin{corollary}  \label{cor-34a}     
	Let $\ch K \ne 2$.                 
	Let $\varphi \in \GL(V)$ be bireflectional. Let $\sigma \in \GL(V)$ be an involution reversing $\varphi$. Assume that $\varphi$ has no elementary devisor $(x \pm 1)^d$ of odd degree $d$. Then $\dim \Fix(\sigma) = \dim \Neg(\sigma) = \frac{\dim V}{2}$.
\end{corollary}

\begin{proof}[Proof of \ref{app-th1}]
Cleary, $\varphi$ is bireflectional if $\varphi$ has an elementary divisor $(x \pm 1)^d$ of odd degree. Otherwise, $\varphi$ is bireflectional if and only if
$\dim V \equiv 0 \mod 4$ by \ref{cor-34a}.
\end{proof}

Using \ref{lemma-inv5}, we obtain the following result of 
Gill and Singh \cite[Proposition 5.5]{GillSingh-2011}:

\begin{theorem} \label{app-th3}
	Let $K = \GF(q)$ be finite. 
	Let $\varphi \in \SL(4m+2, K)$ be similar to its inverse in $\GL(n, K)$. 
	Then the following are equivalent
	\begin{enumerate}
		\item  $q \equiv 1 \mod 4$ or $\varphi$ has an elementary divisor that is a nonsquare;
		\item $\varphi$ is conjugate  to its inverse in $\SL(4m+2, K)$.
	\end{enumerate}
\end{theorem}

%%%%%%%%%%%%%%%%%%%%%%%%%%%%%%%%%%%%%%%%%%%%%%%%%%%%%%%%%%%%%

\end{document}

\typeout{get arXiv to do 4 passes: Label(s) may have changed. Rerun}